\documentclass[a4paper,12pt]{article}

\usepackage{amsthm}
\usepackage{amsmath,amssymb,latexsym,amsfonts,mathrsfs}

\theoremstyle{plain}
\newtheorem{Thm}{Theorem}[section]

\newtheorem{Lem}[Thm]{Lemma}
\newtheorem{Prop}[Thm]{Proposition}

\theoremstyle{definition}

\newcommand{\Proof}[2][Proof]{\begin{proof}[{#1}] #2 \end{proof}}

\numberwithin{equation}{section}

\makeatletter
\renewcommand\section{\@startsection {section}{1}{\z@}%
                                   {-3.5ex \@plus -1ex \@minus -.2ex}%
                                   {2.3ex \@plus.2ex}%
                                   {\normalfont\large\bf}}
\makeatother

\makeatletter
\renewcommand\subsection{\@startsection {subsection}{1}{\z@}%
                                   {-3.5ex \@plus -1ex \@minus -.2ex}%
                                   {2.3ex \@plus.2ex}%
                                   {\normalfont\normalsize\bf}}
\makeatother

\makeatletter
\@addtoreset{footnote}{page}
\makeatother

\newcommand{\rbra}[1]{\left( #1 \right)} 
\newcommand{\cbra}[1]{\left\{ #1 \right\}} 

\renewcommand{\d}{{\rm d}} 

\newcommand{\tend}[2]{\mathrel{\mathop{\longrightarrow}\limits^{#1}_{#2}}}

\renewcommand{\tilde}{\widetilde}

\newcommand{\bN}{\ensuremath{\mathbb{N}}}

\newcommand{\cB}{\ensuremath{\mathcal{B}}}
\newcommand{\cF}{\ensuremath{\mathcal{F}}}

\newcommand{\cP}{\ensuremath{\mathcal{P}}}

\begin{document}
\begin{center}
{\Large \bf Strong solutions of Tsirelson's equation in discrete time 
taking values in compact spaces with semigroup action 
}
\end{center}

\begin{center}
Takao Hirayama\footnote{
The Bank of Tokyo-Mitsubishi UFJ, Ltd, 
2-7-1 Marunouchi, Chiyoda-ku, Tokyo, 100-8388 JAPAN.}
\quad and \quad 
Kouji Yano\footnote{
Graduate School of Science, Kyoto University, Kyoto, JAPAN.}\footnote{
Partially supported by KAKENHI (20740060) and by Inamori Foundation.} 
\end{center}

\begin{abstract}
Under the assumption that the infinite product of evolution process converges almost surely, 
the set of strong solutions are characterized by a compact space, 
which may be regarded as the set of possible initial states. 
\end{abstract}

\noindent
{\footnotesize Keywords and phrases: Stochastic equation, strong solution, 
infinite product of random variables.} 
\\
{\footnotesize AMS 2010 subject classifications: 
Primary
60B15; 
secondary
60J10; 
37H10. 
}

\section{Introduction}

Let $ S $ and $ \Sigma $ be compact metric spaces with countable bases 
and suppose that $ \Sigma $ is a topological semigroup acting continuously on $ S $. 
Denote $ \bN = \{ 0,1,2,\ldots \} $ and consider the following stochastic equation 
(which we call {\em Tsirelson's equation in discrete time}): 
\begin{align}
X_k = N_k X_{k-1} 
\quad \text{for} \ k \in -\bN, 
\label{eq: SE}
\end{align}
where $ X=(X_k)_{k \in -\bN} $ is the unknown process taking values in $ S $ 
and $ N=(N_k)_{k \in -\bN} $ is the driving noise process taking values in $ \Sigma $. 
More precisely, for a given sequence $ (\mu_k)_{k \in -\bN} $ of laws on $ \Sigma $, 
the process 
\begin{align}
\{ (X_k)_{k \in -\bN},(N_k)_{k \in -\bN} \} 
\label{eq: sol}
\end{align}
is called a {\em solution} of equation \eqref{eq: SE} if 
it satisfies \eqref{eq: SE} and for each $ k \in -\bN $ 
the random variable $ N_k $ has law $ \mu_k $ and is independent 
of $ \sigma(X_j:j \le k-1) $. 
We adopt the convention that two solutions are identified 
if their joint laws on $ S^{-\bN} \times \Sigma^{-\bN} $ are equal in law. 
For instance, uniqueness always means the uniqueness in law. 
The process $ (X_k)_{k \in -\bN} $ evolves forward in time $ k $ so that 
the present state $ X_k $ is obtained from $ X_{k-1} $, the state one step before, 
by being acted by $ N_k $. 
Here we must keep in mind that 
the index $ k $ varies in $ -\bN $, the set of non-positive integers, 
so that 
there are a priori no initial time nor initial state in this evolution. 

Equation \eqref{eq: SE} on the one-dimensional torus with trivial group action 
has been originally studied by Tsirelson \cite{MR0375461} 
and further studied by Yor \cite{MR1147613}. 
It was generalized to compact groups by Akahori--Uenishi--Yano \cite{MR2365485} 
and by the authors \cite{HY} (see \cite{YY} for the related survey). 
It was generalized to compact spaces with semigroup action by Yano \cite{Yrcp}. 
For other related works, see Takahashi \cite{MR2582432}, 
Raja \cite{Raja}, Evans--Gordeeva \cite{EV} and Delattre--Rosenbaum \cite{DR}. 
Note also that equation \eqref{eq: SE} can also be found in Furstenberg \cite{F}.

A solution \eqref{eq: sol} of equation \eqref{eq: SE} is called {\em strong} 
if, for each $ k \in -\bN $, 
the present state $ X_k $ is measurable with respect to the past noise $ N_k,N_{k-1},\ldots $ 
up to null sets; 
or in other words, there exists a Borel function $ f_k:\Sigma^{-\bN} \to S $ such that 
\begin{align}
X_k = f_k(N_k,N_{k-1},\ldots) 
\quad \text{almost surely.} 
\label{}
\end{align}

Let $ (\mu_k)_{k \in -\bN} $ be a sequence of laws on $ \Sigma $ 
and let $ (N_k)_{k \in -\bN} $ be a sequence of independent random variables 
such that $ N_k $ has law $ \mu_k $ for each $ k \in -\bN $. 
If the infinite product 
\begin{align}
\Phi_k = \lim_{l \to -\infty } N_k N_{k-1} \cdots N_{l+1} 
\label{eq: infin prod}
\end{align}
converges almost surely for each $ k \in -\bN $, then, for each $ x \in S $, we see that 
\begin{align}
\{ (\Phi_k x)_{k \in -\bN},(N_k)_{k \in -\bN} \} 
\label{eq: str sol}
\end{align}
is a strong solution. 
For this solution, one may think that the point $ x $ might almost be an initial state, 
but, in fact, we cannot do so in certain cases such as Theorem \ref{thm2} as below. 
The following is the main result of this paper, 
which suggests an alternative to the initial states. 

\begin{Thm} \label{thm: main}
Suppose that the infinite product \eqref{eq: infin prod} converges almost surely 
for each $ k \in -\bN $. 
Then there exist a compact Hausdorff space $ T $ with a coutable base, 
a continuous onto mapping $ \pi:S \to T $ 
and a measurable section $ \psi:T \to S $ of $ \pi $ 
(i.e., $ \pi \circ \psi $ is identity) 
which satisfy the following conditions: 
\begin{enumerate}
\item 
for any two distinct elements $ y_1,y_2 \in T $, 
the solutions \eqref{eq: str sol} for $ x=\psi(y_1) $ and $ x=\psi(y_2) $ are distinct; 
\item 
the solutions \eqref{eq: str sol} for $ x=\psi(y) $ with $ y $ running over $ T $ 
exhaust all strong solutions; 
\item 
any solution \eqref{eq: sol} is equal to 
\begin{align}
\{ (\Phi_k \psi(\Xi))_{k \in -\bN},(N_k)_{k \in -\bN} \} 
\label{}
\end{align}
for some $ T $-valued random variable $ \Xi $ which is independent of $ (N_k)_{k \in -\bN} $. 
\end{enumerate}
\end{Thm}

Theorem \ref{thm: main} will be proved in Section \ref{sec: main}. 

Theorem \ref{thm: main} provides us with a general framework 
which unifies the following two earlier studies, which seem in completely different situations. 

$ {\bf 1^{\circ}).} $ 
Suppose that $ S = \Sigma = G $ 
for a compact metric group $ G $ with a countable base. 
(We note that 
the Ellis theorem \cite{MR0088674} asserts that 
a topological semigroup which is algebraically a group 
is necessarily a topological group; 
in particular, the inversion operation is continuous as well.)
We study equation \eqref{eq: SE} 
where $ N_k X_{k-1} $ 
in the right hand side of \eqref{eq: SE} 
is considered to be the usual product in $ G $. 
In \cite{HY}, the authors utilized the results of Csisz\'ar \cite{MR0205306} 
concerning convergence of infinite product of $ G $-valued random variables, 
and obtained the following result: 

\begin{Thm}[\cite{HY}] \label{thm1}
Suppose that there exists a solution \eqref{eq: sol} which is strong. 
Then there exists a sequence $ (\alpha_k)_{k \in -\bN} $ of deterministic elements of $ G $ 
such that the ``centered processes" defined as 
\begin{align}
N^{(\alpha )}_k := \alpha _{k+1}^{-1} N_k \alpha _k 
, \quad 
X^{(\alpha )}_k := \alpha _{k+1}^{-1} X_k 
\quad \text{for} \ k \in -\bN 
\label{}
\end{align}
satisfy the following: 
\begin{enumerate}
\item 
for each $ k \in -\bN $, 
the infinite product $ N^{(\alpha )}_k N^{(\alpha )}_{k-1} \cdots N^{(\alpha )}_{l+1} $ 
converges almost surely as $ l \to -\infty $ to the limit random variable $ \Phi^{(\alpha )}_k $; 
\item 
any strong solution is of the form \eqref{eq: str sol} for some $ x \in G $, 
where $ \Phi_k = \alpha _{k+1} \Phi^{(\alpha )}_k $. 
\end{enumerate}
\end{Thm}

We note that, under the assumptions of Theorem \ref{thm1}, 
Theorem \ref{thm: main} holds with $ T=G $ and $ \pi:G \to G $ being identity. 

$ {\bf 2^{\circ}).} $ 
Suppose that $ S $ consists of finite elements 
and $ \Sigma $ is the composition semigroup of all mappings from $ S $ to itself. 
We equip $ S $ and $ \Sigma $ with discrete topologies. 
We study equation \eqref{eq: SE} 
where $ N_k X_{k-1} $ 
in the right hand side of \eqref{eq: SE} 
is considered to be $ N_k(X_{k-1}) $, the value of the mapping $ N_k $ evaluated at $ X_{k-1} $. 
In \cite{Yrcp}, the second author obtained the following result. 

\begin{Thm}[\cite{Yrcp}] \label{thm2}
Let $ \mu $ be a law on $ \Sigma $. 
Suppose that the Markov chain whose transition probability is 
\begin{align}
p(x,y) = \mu( \sigma : \sigma x = y ) 
\label{}
\end{align}
is ergodic. Set $ \mu_k = \mu $ for all $ k \in -\bN $. 
(In this case, there exists a unique solution.) 
Then the unique solution is strong if and only if 
the {\em road coloring} induced from the support of $ \mu $ is {\em synchronizing}; 
in other words, there exists a finite sequence $ \{ \sigma_0,\sigma_1,\ldots,\sigma_r \} $ 
of the support of $ \mu $ such that 
the composition product $ \sigma_r \sigma_{r-1} \cdots \sigma_0 $ 
maps $ S $ into a singleton. 
In this case, for each $ k \in -\bN $, 
the infinite product \eqref{eq: infin prod} converges almost surely and is given as 
\begin{align}
\Phi_k = N_k N_{k-1} \cdots N_{k-T_k} 
\label{}
\end{align}
where 
\begin{align}
T_k = \inf \cbra{ n \ge r : 
N_{k-n+r} = \sigma_r , \ N_{k-n+r-1} = \sigma_{r-1} , \ldots, \ N_{k-n+0} = \sigma_0 } . 
\label{}
\end{align}
Consequently, the unique strong solution is given as \eqref{eq: str sol} 
for any choice of $ x \in S $. 
\end{Thm}

We note that, under the assumptions of Theorem \ref{thm2}, 
Theorem \ref{thm: main} holds with $ T $ being a singleton. 

Note that Theorem \ref{thm2} is related to the {\em coupling from the past}; 
see, e.g., \cite{H} and also \cite{DF}. 
For other related works, see Yano--Yasutomi \cite{YanoYasu,YanoYasu2}.

This paper is organized as follows. 
In Section \ref{sec: notation}, 
we show that equation \eqref{eq: SE} can be reduced to convolution equation. 
Section \ref{sec: main} is devoted to the proof of Theorem \ref{thm: main}.

\section{Convolution equation} \label{sec: notation}

For general theory of topological semigroups, 
see, e.g., \cite{MR0223483}, \cite{MR0467871} and \cite{MR1363260}. 

Let $ \cB(S) $ denote the set of all Borel sets of $ S $, 
and let $ \cP(S) $ denote the set of all probability laws on $ S $. 
We introduce $ \cB(\Sigma) $ and $ \cP(\Sigma) $ similarly. 
We equip $ \cP(S) $ and $ \cP(\Sigma) $ with the topology of weak convergence. 
Since $ S $ and $ \Sigma $ are compact, they are compactly metrizable. 
For $ \mu_1, \mu_2, \mu \in \cP(\Sigma) $ and $ \lambda \in \cP(S) $, 
we define the {\em convolutions} $ \mu_1 * \mu_2 \in \cP(\Sigma) $ 
and $ \mu * \lambda \in \cP(S) $ by 
\begin{align}
(\mu_1 * \mu_2)(A) =& \int_{\Sigma} \mu_1(\d \sigma_1) \int_{\Sigma} \mu_2(\d \sigma_2) 
1_{\{ \sigma_1 \sigma_2 \in A \}} 
, \quad A \in \cB(\Sigma) , 
\label{} \\
(\mu * \lambda)(A) =& \int_{\Sigma} \mu(\d \sigma) \int_S \lambda(\d x) 
1_{\{ \sigma x \in A \}} 
, \quad A \in \cB(S) . 
\label{}
\end{align}
By the semigroup structure of $ \Sigma $, we see that 
\begin{align}
\mu_1 * (\mu_2 * \mu_3) 
= (\mu_1 * \mu_2) * \mu_3 
, \quad \mu_1,\mu_2,\mu_3 \in \cP(\Sigma) . 
\label{}
\end{align}
By the associativity of the $ \Sigma $-action on $ S $, we see that 
\begin{align}
\mu_1 * (\mu_2 * \lambda) 
= (\mu_1 * \mu_2) * \lambda . 
\label{}
\end{align}
Since $ S $ and $ \Sigma $ are compact, 
the convolutions $ (\mu_1,\mu_2) \mapsto \mu_1 * \mu_2 $ and 
$ (\mu,\lambda) \mapsto \mu * \lambda $ are jointly continuous. 

Let $ (\mu_k)_{k \in -\bN} $ be a sequence of $ \cP(\Sigma) $. 
Let \eqref{eq: sol} be a solution of equation \eqref{eq: SE} 
and let $ (\lambda_k)_{k \in -\bN} $ denote its marginal law system, i.e., 
$ \lambda_k $ is the law of $ X_k $ for each $ k \in -\bN $. 
Then it follows by definition of solution that the convolution equation 
\begin{align}
\lambda_k = \mu_k * \lambda_{k-1} 
\quad \text{for} \ k \in -\bN 
\label{eq: conv eq}
\end{align}
holds. 
The following proposition, which generalizes Lemma 4.3 of \cite{MR2365485}, 
asserts that equation \eqref{eq: SE} can be reduced to 
the convolution equation \eqref{eq: conv eq}. 

\begin{Prop} \label{thm: conv eq}
The following statements hold: 
\begin{enumerate}
\item 
Two solutions whose marginal law systems coincide are equal. 
\item 
Let $ (\lambda_k)_{k \in -\bN} $ be a sequence of $ \cP(S) $ 
such that \eqref{eq: conv eq} holds. Then there uniquely exists a solution \eqref{eq: sol} 
whose marginal law system is $ (\lambda_k)_{k \in -\bN} $. 
\end{enumerate}
\end{Prop}

\Proof{
(i) 
Let two solutions be given such that their marginal law systems coincide. 
Then, by equation \eqref{eq: conv eq}, 
their finite dimensional distributions coincide, 
which implies that they are equal. 

(ii) 
Let $ (\lambda_k)_{k \in -\bN} $ be a sequence of $ \cP(S) $ 
such that \eqref{eq: conv eq} holds. 
For each $ k \in -\bN $, we define a law $ \Lambda_{k+1} $ 
on $ (\Sigma \times S)^{-k} $ as follows. 
Let $ N_0,N_{-1},\ldots,N_{k+1},X_k $ be independent random variables 
such that $ N_j $ has law $ \mu_j $ for $ j=0,-1,\ldots,k+1 $ 
and that $ X_k $ has law $ \lambda_k $. 
For $ j=0,-1,\ldots,k+1 $, we set $ X_j = N_{j,k} X_k $, 
where 
\begin{align}
N_{j,k} = N_j N_{j-1} \cdots N_{k+1} 
\quad \text{for} \ j>k . 
\label{}
\end{align}
Define $ \Lambda_{k+1} $ be the law of $ (X_j,N_j:j=0,-1,\ldots,k+1) $ 
on $ (\Sigma \times S)^{-k} $. 
Then, by \eqref{eq: conv eq}, 
it is obvious that the family $ \{ \Lambda_{k+1}:k \in -\bN \} $ is consistent. 
Thus, by Kolmogorov's extension theorem, we obtain existence of a solution. 
}

\section{Proof of main theorem} \label{sec: main}

Since $ \cP(S) $ is a compact Hausdorff space with a countable base, 
so is the countable direct product space $ \cP(S)^{-\bN} $. 
Since $ \cP(S) $ is equipped with a convex structure in the usual way, 
so is $ \cP(S)^{-\bN} $. 
Let us write $ \cP^{\rm cvl} $ for the set of all solutions $ (\lambda_k)_{k \in -\bN} $ 
of equation \eqref{eq: conv eq}. 

\begin{Lem} \label{lem: cvl}
The set $ \cP^{\rm cvl} $ is a compact convex subset of $ \cP(S)^{-\bN} $. 
\end{Lem}

\Proof{
It is obvious by equation \eqref{eq: conv eq} that 
$ \cP^{\rm cvl} $ is a convex subset of $ \cP(S)^{-\bN} $. 
Let us prove that $ \cP^{\rm cvl} $ is closed in $ \cP(S^{-\bN}) $. 
Let $ (\lambda^{(n)}_k)_{k \in -\bN} \subset \cP^{\rm cvl} $ 
be such that $ (\lambda^{(n)}_k)_{k \in -\bN} \to (\lambda_k)_{k \in -\bN} $ 
as $ n \to \infty $ for some $ (\lambda_k)_{k \in -\bN} \in \cP(S)^{-\bN} $. 
Then we have 
\begin{align}
\lambda^{(n)}_k = \mu_k * \lambda^{(n)}_{k-1} 
\quad \text{for all} \ k \in -\bN \ \text{and} \ n \in \bN . 
\label{}
\end{align}
Since $ \lambda^{(n)}_k \to \lambda_k $ as $ n \to \infty $ for all $ k \in -\bN $, 
we see, by the continuity of convolutions, that equation \eqref{eq: conv eq} holds, 
which implies that $ (\lambda_k)_{k \in -\bN} \in \cP^{\rm cvl} $. 
Since $ \cP(S)^{-\bN} $ has a countable base, 
we conclude that $ \cP^{\rm cvl} $ is closed. 
}

Now we proceed to prove Theorem \ref{thm: main}. 

\Proof[Proof of Theorem \ref{thm: main}]{
Let us assume that 
the infinite product \eqref{eq: infin prod} converges almost surely 
for each $ k \in -\bN $. 
We then see that 
\begin{align}
\mu_{k,l} := \mu_k * \mu_{k-1} * \cdots * \mu_{l+1} 
\tend{}{l \to -\infty } 
\nu_k 
\quad \text{for all} \ k \in -\bN , 
\label{eq: repre1}
\end{align}
where $ \nu_k $ is the law of $ \Phi_k $ for all $ k \in -\bN $. 

Let $ T = \cP^{\rm cvl} $. By Lemma \ref{lem: cvl}, 
we see that $ T $ is a compact Hausdorff space with a countable base. 
For $ x \in S $, we write $ \pi(x) = (\nu_k * \delta_x)_{k \in -\bN} $. 
It is now obvious that the mapping $ \pi:S \to T $ is continuous. 
By Theorem 6.9.7 of \cite{B}, there exists a measurable section $ \psi:T \to S $ of $ \pi $. 
Claim (i) is obvious by Proposition \ref{thm: conv eq}.

Let $ \{ (X_k)_{k \in -\bN},(N_k)_{k \in -\bN} \} $ be an arbitrary strong solution. 
Let $ \lambda_k $ denote the law of $ X_k $ for all $ k \in -\bN $. 
We then have $ (\lambda_k)_{k \in -\bN} \in \cP^{\rm cvl} $. 
Since we have 
\begin{align}
X_k = N_{k,l} X_l 
\label{eq: eqeq2}
\end{align}
by iterating equation \eqref{eq: SE}, 
the conditional law of $ X_k $ given $ \cF_{k,l} := \sigma(N_k,N_{k-1},\cdots,N_{l+1}) $ 
is written as 
\begin{align}
P(X_k \in \cdot|\cF_{k,l}) = \delta_{N_{k,l}} * \lambda_l . 
\label{}
\end{align}
Letting $ l \to -\infty $, we see that 
the left hand side converges to $ \delta_{X_k} $ 
since $ X_k $ is measurable with respect to $ \cF_{k,-\infty } := \sigma(N_k,N_{k-1},\ldots) $ 
up to null sets. 
Taking a subsequence if necessary, we may assume that $ \lambda_l \to \lambda $ 
for some $ \lambda \in \cP(S) $, 
we obtain 
\begin{align}
\delta_{X_k} = \delta_{\Phi_k} * \lambda 
= \int_S \lambda(\d x) \rbra{ \delta_{\Phi_k} * \delta_x } . 
\label{}
\end{align}
This shows that 
\begin{align}
(\lambda_k)_{k \in -\bN} = \pi(x) 
\quad \text{for $ \lambda $-almost every $ x \in S $}, 
\label{}
\end{align}
which implies, in particular, that 
$ (\lambda_k)_{k \in -\bN} \in T $. 
This proves Claim (ii).

Let $ \{ (X_k)_{k \in -\bN},(N_k)_{k \in -\bN} \} $ be an arbitrary solution. 
Let $ \lambda_k $ denote the law of $ X_k $ for all $ k \in -\bN $. 
We then have $ (\lambda_k)_{k \in -\bN} \in \cP^{\rm cvl} $. 
Since we have \eqref{eq: eqeq2} by iterating equation \eqref{eq: SE}, 
we have 
\begin{align}
\lambda_k = \mu_{k,l} * \lambda_l 
\quad \text{for} \ k>l . 
\label{}
\end{align}
Taking a subsequence if necessary, we may assume that $ \lambda_l \to \lambda $ 
for some $ \lambda \in \cP(S) $, so that we obtain 
\begin{align}
\lambda_k = \nu_k * \lambda 
= \int_S \lambda(\d x) (\nu_k * \delta_x) 
\quad \text{for} \ k \in -\bN . 
\label{}
\end{align}
Now we obtain 
\begin{align}
(\lambda_k)_{k \in -\bN} = \int_S \lambda(\d x) \pi(x) = \int_T \tilde{\lambda}(\d y) y , 
\label{}
\end{align}
where $ \tilde{\lambda} = \lambda \circ \pi^{-1} $. 
Taking an independent random variable $ \Xi $ whose law is $ \tilde{\lambda} $, 
we obtain Claim (iii). 

The proof is now complete. 
}

{\bf Acknowledgements:} 
The first author, T.~ H., would like to express his sincere thanks 
to Professor Jir\^o Akahori, who kindly expended a considerable effort 
to support T.~H. in his study. 

\def\cprime{$'$} \def\cprime{$'$}

\end{document}